\documentclass[a4paper,10pt,reqno]{amsart}
\date{November 9, 2026}
\usepackage{graphicx}
\usepackage{latexsym}
\usepackage{amssymb}
\usepackage{amsmath}
\usepackage{amscd}
\usepackage{amsthm}
\usepackage{amsfonts}
\usepackage{mathrsfs}
\usepackage{enumerate}
\usepackage{enumitem}
\usepackage{bm}
\usepackage[usenames]{color}
\usepackage[active]{srcltx}
\usepackage{comment}
\usepackage{here}
\usepackage[normalem]{ulem}
\usepackage[subrefformat=parens]{subcaption}
  \DeclareSymbolFontAlphabet{\mathbb}{AMSb}
\title[Curve theory in $S^2\times S^2$]%
{Curve theory in $S^2\times S^2$
with applications to ruled surfaces in $S^3$}
\author[A.~Honda]{Atsufumi Honda}
\thanks{Corresponding author: Atsufumi Honda.}
\address{%
   Department of Applied Mathematics, 
   Faculty of Engineering, Yokohama National University, %\endgraf
   Hodogaya, Yokohama 240-8501, Japan
}
\email{honda-atsufumi-kp@ynu.ac.jp}
\author[J.~Inoguchi]{Jun-ichi Inoguchi}
\address{%
   Department of Mathematics, 
   Hokkaido University, 
   Sapporo, 060-0810, Japan}
\email{inoguchi@math.sci.hokudai.ac.jp}
\author[H.~Kiyohara]{Hirotaka Kiyohara}
\address{%
   Department of Mathematics Education, 
   Osaka Kyoiku University, 
   Osaka, 582-8582, Japan}
\email{kiyohara-h51@cc.osaka-kyoiku.ac.jp}
\dedicatory{Dedicated to Professor Udo Hertrich-Jeromin and 
Professor Wayne Rossman
on the occasions of their 60th birthdays}
\theoremstyle{plain}
 \newtheorem{theorem}{Theorem}[section]
 \newtheorem{introtheorem}{Theorem}
 \newtheorem{introcorollary}[introtheorem]{Corollary}

 \newtheorem{proposition}[theorem]{Proposition}
 
 \newtheorem*{fact*}{Fact}
 \newtheorem{lemma}[theorem]{Lemma}
 \newtheorem{corollary}[theorem]{Corollary}
 \theoremstyle{remark}
 \newtheorem{definition}[theorem]{Definition}
 \newtheorem{remark}[theorem]{Remark}
 
 \newtheorem{example}[theorem]{Example}
\numberwithin{equation}{section}

\newcommand{\N}{\mathbb{N}}%{\boldsymbol{N}}%
\newcommand{\Z}{\mathbb{Z}}%{\boldsymbol{Z}}%
\newcommand{\Q}{\mathbb{Q}}
\newcommand{\R}{\mathbb{R}}%{\R}%
\newcommand{\C}{\mathbb{C}}%{\C}%
\renewcommand{\H}{\mathbb{H}}%{\boldsymbol{H}}
\renewcommand{\S}{S}%{\boldsymbol{S}}

\newcommand{\SU}{\operatorname{SU}}
\newcommand{\SO}{\operatorname{SO}}
\newcommand{\Isom}{\operatorname{Isom}}

\newcommand{\vect}[1]{\boldsymbol{#1}}

\newcommand{\LS}{\mathcal{L}(\S^3)}

\renewcommand{\Im}{\operatorname{Im}}

\newcommand{\e}{\boldsymbol{e}}
\newcommand{\n}{\boldsymbol{n}}
\newcommand{\TC}{\operatorname{TC}}

\begin{document}

\begin{abstract}
The Riemannian product $S^2\times S^2$ of two unit spheres
admits several important geometric interpretations.
It can be identified both with the complex quadric $Q^2$
and with the space $\LS$ of oriented geodesics
in the unit $3$-sphere $S^3$.
In this paper, we first develop a unified framework
for the geometry of curves in $S^2\times S^2$.
We introduce projection-regular curves in $S^2\times S^2$
and their curvature data,
and prove a congruence theorem for these curves.
As an application of this framework,
regarding $S^2\times S^2$ as the complex quadric $Q^2$,
we classify K\"ahler magnetic geodesics
associated with the product K\"ahler--Einstein structure.
For closed generic projection-regular curves,
we establish a sharp Fenchel-type inequality.
We then identify $S^2\times S^2$ with the space $\LS$
of oriented geodesics in the unit $3$-sphere $S^3$
and apply our curve theory to ruled surfaces in $S^3$.
In particular, we prove that every nonsingular closed
generic projection-regular ruled surface attaining equality
in the Fenchel-type inequality is an embedded torus.
\end{abstract}
\maketitle

%\tableofcontents

\section{Introduction}
\label{sec:introduction}

This paper develops a curve theory 
adapted to the product geometry of
$S^2\times S^2$.
Here, $S^2$ denotes the (Euclidean) unit $2$-sphere.
A curve
\[
 \gamma=(\gamma_1,\gamma_2)\colon I\longrightarrow S^2\times S^2
\]
has two spherical projections, and their relative speeds
and signed curvatures contain information that is invisible to the
ordinary Frenet curvatures of the curve in 
$S^2\times S^2$.
The first part of the paper formulates this information.
The second part identifies $S^2\times S^2$ with the
space of oriented geodesics in 
the unit $3$-sphere $S^3$ and applies 
the results of the first part to ruled surfaces.

\subsection{Projection-regular curves in $S^2\times S^2$}
A unit-speed Frenet curve in a simply connected real space form 
is uniquely determined, up to congruence, by its Frenet curvatures.
In a general Riemannian manifold, however, the Frenet curvatures
are \emph{not} sufficient, and additional invariants adapted to the
geometry of the ambient manifold are required
(see, for example, \cite{CM,CMM}).
For instance, Maeda--Ohnita \cite{Maeda-Ohnita} showed that 
curves in complex space forms are determined by their 
Frenet curvatures and complex torsions.
Kimura--Ortega \cite{Kimura-Ortega} established a congruence
theorem for Frenet curves of osculating order two
in the complex quadric
$Q^2\simeq \S^2\times\S^2$
(cf.~Remark~\ref{rem:frenet-triplet}).

In this paper, we take a different approach.
We introduce the class of \emph{projection-regular curves}.  
A regular curve in $S^2\times S^2$ is 
called \emph{projection-regular}
if each of its two projections is either regular or constant.  
It is \emph{generic}
(G-projection-regular) when both projections are regular 
and \emph{exceptional} (E-projection-regular)
when exactly one projection is regular.  

For a unit-speed G-projection-regular curve, 
the component speeds can be written uniquely
as $\cos\theta$ and $\sin\theta$, where
$0<\theta<\pi/2$.  Together with the signed curvature coefficients
$\kappa_1$ and $\kappa_2$ of the two projected curves, this gives the
\emph{curvature triplet}
$(\theta,\kappa_1,\kappa_2)$.
For an E-projection-regular curve, the curvature data consist of its type and the signed
curvature coefficient $\kappa_0$ of its regular projection.  Our first
main result shows that these data form a complete system of congruence
invariants for unit-speed projection-regular curves.

\begin{introtheorem}[Congruence theorem]\label{thm:A}
For any smooth functions $\theta,\kappa_1,\kappa_2$ on an interval, with
$0<\theta<\pi/2$, there exists a unit-speed G-projection-regular curve in
$S^2\times S^2$ having
curvature triplet $(\theta,\kappa_1,\kappa_2)$, and it is unique up to
congruence.  The analogous existence and uniqueness statement holds for
E-projection-regular curves of either type with prescribed curvature $\kappa_0$.
\end{introtheorem}

A more detailed statement of Theorem~\ref{thm:A} is given in
Theorem~\ref{thm:fundamental-thm}.
Note that a G-projection-regular curve may have
osculating order three or four
(see Remark~\ref{rem:frenet-triplet}).
Consequently, our congruence theorem is not a special case of the
congruence theorem in \cite{Kimura-Ortega}.

As an application of our theory of projection-regular curves,
we obtain a classification of 
K\"ahler magnetic geodesics 
in $\S^2\times\S^2\simeq Q^2$
with respect to the canonical K\"ahler--Einstein structure
(Theorem~\ref{thm:Kahler-magnetic}).

\subsection{Total curvature estimates for closed curves in $S^2\times S^2$}
For a regular closed curve $\gamma\colon S^1\to S^2$ 
in the unit $2$-sphere $S^2$,
its total absolute geodesic curvature
$\int_{S^1}|k_g|\,ds$ is nonnegative by definition, and equality is
attained by geodesics.  In this sense, the resulting inequality is
trivial in comparison with its Euclidean counterpart.
On the other hand, regarding $\gamma$ as a closed curve in 
Euclidean $3$-space $\R^3$,
its total absolute curvature is
$\int_{S^1}\sqrt{1+k_g^2}\,ds$.  
By Fenchel's theorem \cite{Fenchel1929}, 
this quantity is at least $2\pi$, 
with equality if and only if
$\gamma\colon S^1\to S^2$ is an embedding whose image is a circle.
From the viewpoint of spherical curve geometry, consider the
canonical lift of $\gamma$ to $US^2$
$$
\hat\gamma:=(\gamma,\vect{e}) \colon S^1 \to US^2\qquad
(\vect{e}:=\gamma'/\|\gamma'\|).
$$
Its length with respect to the Sasaki metric coincides with
$\int_{S^1}\sqrt{1+k_g^2}\,ds.$

Likewise, for a curve in $S^2\times S^2$, the integral of its first
Frenet curvature is nonnegative, and equality is attained by
geodesics.  Motivated by the preceding observation, we define the
total curvature of a closed G-projection-regular curve
$\gamma=(\gamma_1,\gamma_2)\colon S^1\to S^2\times S^2$ as follows.
Consider its canonical lift to $US^2\times US^2$
$$
\hat\gamma:=((\gamma_1,\vect{e}_1),(\gamma_2,\vect{e}_2)) \colon
S^1 \to US^2 \times US^2\qquad
(\vect{e}_j:=\gamma_j'/\|\gamma_j'\|, ~j=1,2)
$$
and its length with respect to the product Sasaki metric,
$$
 \mathcal T(\gamma)
 :=\int_{S^1} \|\hat\gamma'\|_{\mathrm{Sas}}\,ds
$$
which we call the \emph{total curvature} of $\gamma$.
We then obtain the following theorem.

\begin{introtheorem}[Fenchel-type inequality]\label{thm:C}
Let $\gamma\colon  S^1\to S^2\times S^2$ be a closed G-projection-regular curve.
Then
\[
 \mathcal T(\gamma)\geq2\sqrt2\,\pi.
\]
Equality holds if and only if 
both projected curves $\gamma_j\colon S^1\to S^2$ $(j=1,2)$ are embeddings
and the curvature triplet is constant.
In particular, in the equality case, 
both images $\gamma_1(S^1)$ and $\gamma_2(S^1)$ are
circles in $S^2$.
\end{introtheorem}

\subsection{Applications to ruled surfaces in $S^3$}
In the final part of this paper, we use the equivariant identification
\[
\LS\cong S^2\times S^2,
\]
where $\LS$ is the space of oriented geodesics in $S^3$.
A one-parameter family of oriented geodesics in $S^3$ is thereby a curve
in $S^2\times S^2$, and hence determines a ruled surface in $S^3$
(Definition~\ref{def:ruled-surface}).
Classical studies of ruled surfaces in $S^3$ include Lawson helicoids
\cite{Lawson}, Hopf tori \cite{Pinkall}, and 
extrinsically flat ruled surfaces \cite{IzumiyaNagaiSaji,Honda2017}.
Our approach differs in that it begins with invariants of 
the one-parameter family of oriented geodesics in $S^3$ itself and then
translates them into local and global properties of the corresponding
ruled surface.

Ruled surfaces corresponding to curves with constant
curvature data in $S^2\times S^2$ are precisely the helicoidal ruled
surfaces in $S^3$, including ruled surfaces of revolution
(Proposition~\ref{prop:helicoidal-characterization}).

For a ruled surface in $S^3$, we define 
the \emph{left and right ruling Gauss maps} 
(see Definition~\ref{def:ruling-Gaussmap}).
Up to a fixed reflection, the pair of ruling Gauss maps agrees with
the curve $\gamma=(\gamma_1,\gamma_2)$ in $S^2\times S^2$
corresponding to the ruled surface
(Proposition~\ref{prop:ruling-Gauss-curve}).
A ruled surface is called \emph{G-projection-regular}
if the associated curve $\gamma$ is G-projection-regular
(Definition~\ref{def:projection-regular-ruled-surface}).
For a closed G-projection-regular ruled surface $F$, 
we define its ruling total curvature
$\mathcal T(F)$ to be the total curvature of $\gamma$.
Theorem~\ref{thm:C} then yields the following corollary.

\begin{introcorollary}\label{cor:C}
Let $F\colon S^1\times S^1\to S^3$ be a closed G-projection-regular ruled surface.
Then
\[
 \mathcal T(F)\geq2\sqrt2\,\pi.
\]
Equality holds if and only if $F$ is simple-rotational.
Moreover, if equality holds and $F$ has no singular points, then
$F$ is an embedding.
\end{introcorollary}

Here, a ruled surface in $S^3$ is called 
\emph{simple-rotational}
if its oriented rulings are generated by a rotational motion 
and the given period of the ruling parameter 
is the least positive period of the generating motion
(see Definition~\ref{def:simple-rotational}).

\begin{figure}[htb]
\begin{center}
 \begin{tabular}{{c@{\hspace{8mm}}c}}
  \resizebox{3.5cm}{!}{\includegraphics{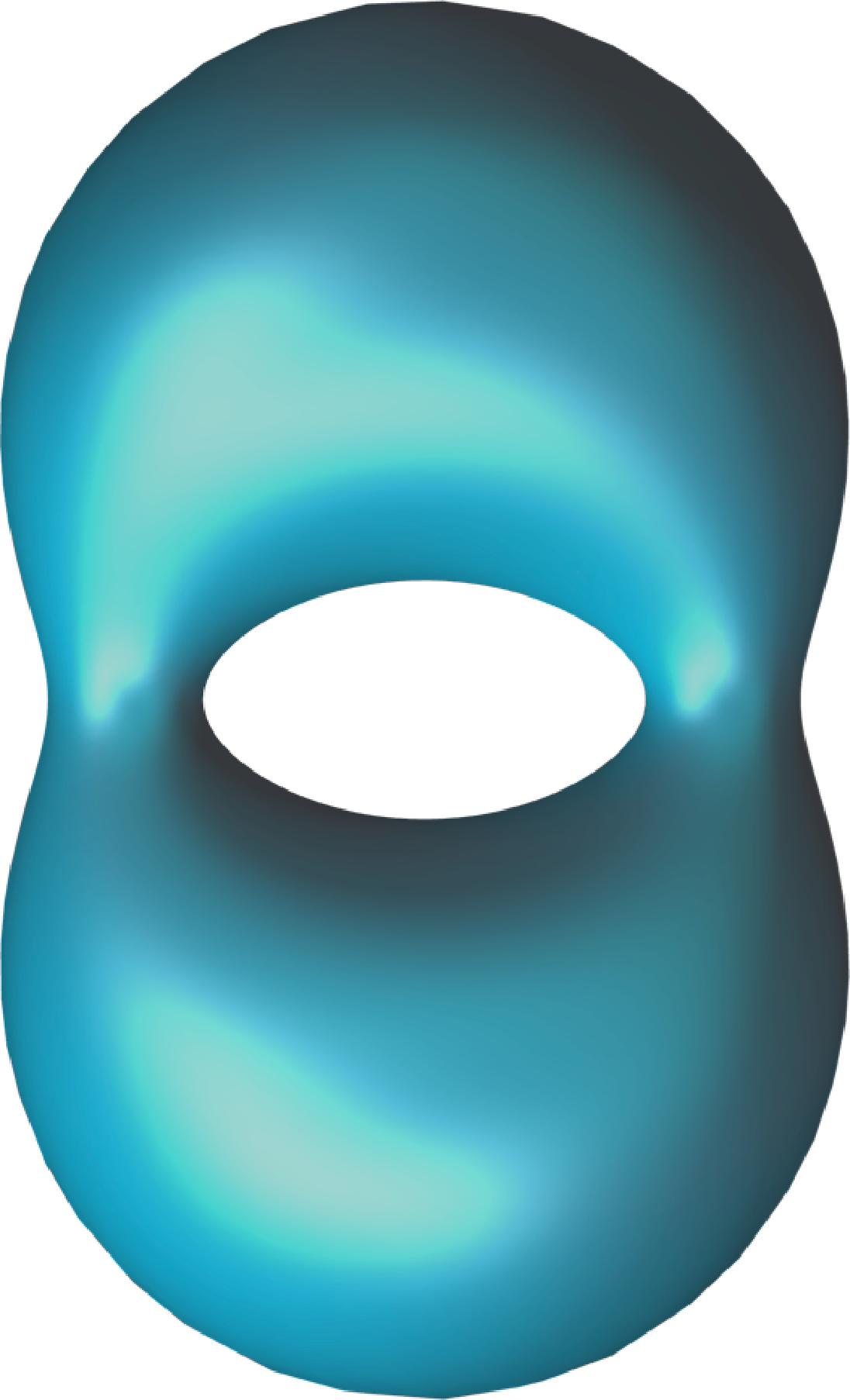}}&
  \resizebox{3.5cm}{!}{\includegraphics{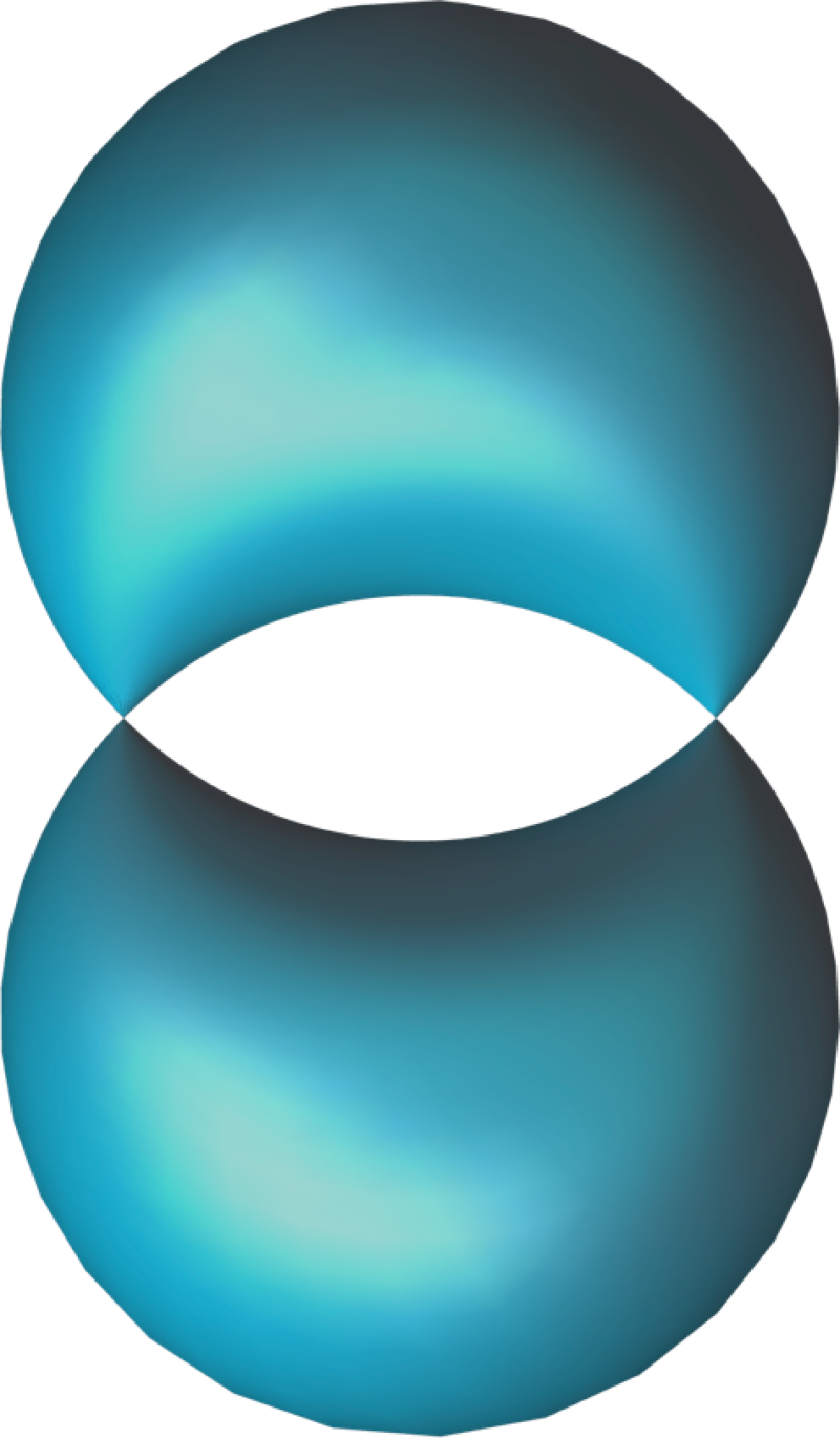}}
 \end{tabular}
\caption{G-projection-regular ruled surfaces in $\S^3$ attaining equality
in the inequality of Corollary~\ref{cor:C}.
The surface on the left is nonsingular and embedded.
The surface on the right has conelike singularities.
Its singular set is homeomorphic to $S^1\sqcup S^1$,
and each of its two components is collapsed to a single point.
In particular, the ruled surface on the right is not injective.}
\label{fig:equality-case}
\end{center}
\end{figure}

\subsection{Organization of this paper}
This paper is organized as follows.
Section~\ref{sec:congruence-thm} introduces projection-regular curves and their curvature
data, proves the congruence theorem, compares the resulting theory with
Frenet theory on $Q^2$, characterizes homogeneous curves, and classifies
K\"ahler magnetic geodesics for the product K\"ahler--Einstein
structure.
Section~\ref{sec:fenchel-type} proves the Fenchel-type inequality for
closed G-projection-regular curves.
Section~\ref{sec:LS3} introduces the identification with $\LS$ and
transfers the preceding results to ruled surfaces in $S^3$.

\section{Projection-regular curves in $S^2\times S^2$}
\label{sec:congruence-thm}

This section develops the basic curve theory used throughout the paper.
We first define projection-regular curves and the curvature data determined by their
spherical projections.  We next prove the corresponding congruence
theorem and compare curves with constant curvature data with Frenet
curves in the complex quadric.  We then characterize projection-regular curves with
constant curvature data as homogeneous curves and give their normal
forms.  Finally, we classify K\"ahler magnetic geodesics using the
preceding curve theory.

\subsection{Projection-regular curves and curvature data}
\label{subsec:projection-regular-data}

In this subsection, we define projection-regular curves and
introduce the curvature data adapted to the two projections.

We regard $S^2$ as the unit sphere in Euclidean space $\R^3$
and equip $M:=S^2\times S^2$ with the product metric
\[
 g=g_{S^2}\oplus g_{S^2}.
\]
The identity component of the isometry group is
\[
 K:=\Isom_0(M)=\SO(3)\times\SO(3).
\]
It acts factorwise on $M$.

For a regular spherical curve $\alpha\colon I\to S^2$, set
$$
v=\|\alpha'\|>0,\qquad
\kappa=\langle\e',\n\rangle,
\qquad
\e=\frac{\alpha'}{v},\qquad
\n=\alpha\times\e.
$$
Then
\begin{equation}\label{eq:spherical-frame}
\alpha'=v\e,\qquad
\e'=-v\alpha+\kappa\n,\qquad
\n'=-\kappa\e
\end{equation}
hold.
The functions $v$ and $\kappa$ are invariant under the action of $\SO(3)$ given by
$$
\alpha\longmapsto A\alpha
\qquad (A\in\SO(3)).
$$
Conversely, the following holds.

\begin{lemma}\label{lem:spherical-frame}
Let $I$ be an open interval.
For given smooth functions $v, \kappa\colon I\to\R$ satisfying $v>0$,
there is a unique regular curve $\alpha\colon I\to S^2$ up to congruence
such that
$v= \|\alpha'\|$ and $\kappa=\langle\e',\n\rangle$.
\end{lemma}

\begin{proof}
For an initial matrix $Q(s_0)\in\SO(3)$, solve
\[
 Q'=Q
 \begin{pmatrix}
 0&-v&0\\
 v&0&-\kappa\\
 0&\kappa&0
 \end{pmatrix}.
\]
The coefficient matrix is skew-symmetric, so $Q(s)$ remains in
$\SO(3)$.  Its three columns satisfy \eqref{eq:spherical-frame}
and the
first column is the required curve.
\end{proof}

\begin{definition}\label{def:projection-regular-curve}
Let $\gamma=(\gamma_1,\gamma_2)\colon I\to S^2\times S^2$ be a regular curve.
It is called a \emph{projection-regular curve} if
each $\gamma_j$ is either regular or constant.
\begin{enumerate}
\item It is \emph{generic} (G-projection-regular) if both projections are regular.
\item It is \emph{exceptional of type $1$} if $\gamma_1$ is regular and
      $\gamma_2$ is constant.
\item It is \emph{exceptional of type $2$} if $\gamma_1$ is constant and
      $\gamma_2$ is regular.
\end{enumerate}
A curve of either exceptional type 
is called \emph{E-projection-regular}.
\end{definition}

Let $\gamma$ be a G-projection-regular curve parametrized by arc length.  Put
\[
 \e_j=\frac{\gamma_j'}{\|\gamma_j'\|},
 \qquad
 \n_j=\gamma_j\times\e_j.
\]
There is a unique smooth function $\theta\colon I\to(0,\pi/2)$ such that
\[
 \|\gamma_1'\|=\cos\theta,\qquad
 \|\gamma_2'\|=\sin\theta.
\]
Define $\kappa_j=\langle\e_j',\n_j\rangle$. 
The frames satisfy the following 
Frenet--Serret-type equations:
\begin{align}\label{eq:FS}
\begin{split}
 \gamma_1'&=\cos\theta\,\e_1,\qquad
 \e_1'=-\cos\theta\,\gamma_1+\kappa_1\n_1,\qquad
 \n_1'=-\kappa_1\e_1,\\
 \gamma_2'&=\sin\theta\,\e_2,\qquad
 \e_2'=-\sin\theta\,\gamma_2+\kappa_2\n_2,\qquad
 \n_2'=-\kappa_2\e_2.
\end{split}
\end{align}
We call $(\theta,\kappa_1,\kappa_2)$ the \emph{curvature triplet}.
Notice that the geodesic curvature of $\gamma_j$ 
as a curve in $S^2$ is $\kappa_j/\|\gamma_j'\|$.

For an E-projection-regular curve, parametrize its regular projection $\gamma_0$ by arc
length and set
\[
 \e_0=\gamma_0',\qquad
 \n_0=\gamma_0\times\e_0,\qquad
 \kappa_0=\langle\e_0',\n_0\rangle.
\]
Its curvature data consist of its type and the function $\kappa_0$.

The curvature data are called \emph{constant} if
$\theta$, $\kappa_1$, and $\kappa_2$ are constant
in the G-projection-regular case, or if $\kappa_0$ is constant
in the E-projection-regular case.

\subsection{The congruence theorem}
\label{subsec:congruence-theorem}

We now prove existence and uniqueness
of projection-regular curves.
Two curves in $M=S^2\times S^2$ are \emph{congruent} if they differ by
one fixed element of $K=\SO(3)\times\SO(3)$.  Equivalently, their two projections differ
by fixed elements of $\SO(3)$.

\begin{theorem}\label{thm:fundamental-thm}
Let $I$ be a nonempty open interval.
\begin{enumerate}
\item Given smooth functions $\theta,\kappa_1,\kappa_2$ on $I$ with
      $0<\theta<\pi/2$, there exists a unit-speed G-projection-regular curve
      whose curvature triplet is $(\theta,\kappa_1,\kappa_2)$.  
      It is
      unique up to congruence.
\item Given a smooth function $\kappa_0$ and a type $r\in\{1,2\}$, there
      exists a unit-speed E-projection-regular curve of type $r$ with curvature
      $\kappa_0$.  It is unique up to congruence among curves of that type.
\end{enumerate}
\end{theorem}

\begin{proof}
Apply Lemma~\ref{lem:spherical-frame} with
$(v,\kappa)=(\cos\theta,\kappa_1)$ and
$(\sin\theta,\kappa_2)$. 
This gives a unique G-projection-regular curve satisfying
\eqref{eq:FS} up to congruence. 
The E-projection-regular case is the same construction 
with one constant component. 
\end{proof}

Note that the two E-projection-regular types 
are not interchanged by $K$,
since its action preserves the two factors.

\begin{remark}[Osculating order of curves with constant curvature data]
\label{rem:frenet-triplet}
For a unit-speed curve $\gamma$ in $S^2\times S^2$, its osculating
order at a point is the dimension of
\[
 \operatorname{span}
 \left\{
 \gamma',
 \nabla_s\gamma',
 \nabla_s^2\gamma',
 \nabla_s^3\gamma',
 \ldots
 \right\},
\]
where $\nabla$ denotes the Levi-Civita connection on $S^2 \times S^2$.
Let $\gamma$ be a G-projection-regular curve with constant curvature triplet
$(\theta,\kappa_1,\kappa_2)$
$(0<\theta<\pi/2)$.
Set
$v_1=\cos\theta$,
$v_2=\sin\theta$.
Using \eqref{eq:FS} and the product Levi--Civita connection,
we have
\begin{gather*}
\gamma' = \bigl(v_1\e_1,v_2\e_2\bigr),\qquad
 \nabla_s\gamma'
 =
 \bigl(v_1\kappa_1\n_1,
       v_2\kappa_2\n_2\bigr),\\
 \nabla_s^2\gamma'
 =
 \bigl(-v_1\kappa_1^2\e_1,
       -v_2\kappa_2^2\e_2\bigr),\qquad
 \nabla_s^3\gamma'
 =
 \bigl(-v_1\kappa_1^3\n_1,
       -v_2\kappa_2^3\n_2\bigr).
\end{gather*}
Consequently, the osculating order is classified as follows:
\begin{enumerate}
 \item
 If $\kappa_1=\kappa_2=0$, 
 then $\gamma$ has osculating order one, 
 that is, $\gamma$ is a geodesic.
\item
If
$
 |\kappa_1|=|\kappa_2|>0,
$
then $\gamma$ has osculating order two.
\item
If exactly one of $\kappa_1$ and $\kappa_2$ vanishes,
then $\gamma$ has osculating order three.
\item
If
$
 \kappa_1\kappa_2
 (|\kappa_1|-|\kappa_2|)\ne0,
$
then $\gamma$ has osculating order four.
\end{enumerate}
Thus, even among curves with constant curvature data, 
the curves covered by the congruence theorem are not, 
in general, order-two Frenet curves. 
Consequently, the present congruence theorem 
(Theorem~\ref{thm:fundamental-thm}) is not
a special case of the congruence theorem for order-two Frenet curves
in \cite{Kimura-Ortega}.  
\end{remark}

\subsection{Homogeneous curves}
\label{sec:homogeneous}

We next characterize constant curvature data in terms of symmetry.
We first prove that 
a regular curve in $S^2\times S^2$ 
parametrized by arc length
is homogeneous if and only if it is projection-regular
with constant curvature data.
We then give an explicit normal form and
describe the closedness and covering numbers of a 
constant-data G-projection-regular curve.

\begin{definition}\label{def:homogeneous-curve}
A curve $\gamma\colon \R\to M$ 
in a Riemannian manifold $M$
is \emph{homogeneous} if
there exists a one-parameter subgroup 
$\{\Phi_s\}\subset K=\Isom_0(M)$
such that
\[
 \gamma(s)=\Phi_s\cdot\gamma(0).
\]
\end{definition}

A regular spherical curve 
$\alpha\colon \R\to S^2$ is homogeneous
if and only if it has constant speed and constant signed
geodesic curvature.
This follows immediately from the constant-coefficient version
of the Frenet--Serret-type 
equations \eqref{eq:spherical-frame}.

\begin{proposition}\label{prop:projection-regular-homogeneous}
Let $\gamma\colon \R\to S^2\times S^2$ be a regular curve parametrized by
arc length.  Then $\gamma$ is a projection-regular curve with constant
curvature data if and only if it is homogeneous.
\end{proposition}

\begin{proof}
If the curvature data are constant, each nonconstant projection has
constant speed and curvature coefficient, hence is an orbit of a
one-parameter subgroup of $\SO(3)$.
The two subgroups together define
a one-parameter subgroup of $K=\SO(3)\times\SO(3)$,
so $\gamma$ is homogeneous. Conversely,
the two projections of a homogeneous curve are rotation orbits.  Each is
either constant or a spherical circle of constant speed and curvature,
which gives constant curvature data.
\end{proof}

\begin{remark}
Related results on orbits of one-parameter subgroups 
in the complex quadric $Q^2$ were obtained by 
Kimura--Ortega \cite{Kimura-Ortega}.
\end{remark}

\begin{proposition}\label{prop:constant-data-normal-form}
Let $0<\theta<\pi/2$ and $\kappa_1,\kappa_2\in\R$ be constants, and put
\begin{equation}\label{eq:constant-data-frequencies}
 \mu_1=\sqrt{\cos^2\theta+\kappa_1^2},\qquad
 \mu_2=\sqrt{\sin^2\theta+\kappa_2^2}.
\end{equation}
Up to congruence, the G-projection-regular curve with these data is
\begin{align}
\label{eq:circle-1}
\gamma_1(s)
&=\frac1{\mu_1}
 \bigl(\cos\theta\cos(\mu_1s),
       \cos\theta\sin(\mu_1s),\kappa_1\bigr),\\
\label{eq:circle-2}
\gamma_2(s)
&=\frac1{\mu_2}
 \bigl(\sin\theta\cos(\mu_2s),
       \sin\theta\sin(\mu_2s),\kappa_2\bigr).
\end{align}
It is closed if and only if $\mu_1/\mu_2\in\Q$.  If $L$ is its least
positive period, then
\begin{equation}\label{eq:constant-data-period}
 \mu_1L=2\pi m_1,\qquad
 \mu_2L=2\pi m_2
\end{equation}
for relatively prime $m_1,m_2\in\N$.  
Let $C_j=\gamma_j(\R)\subset S^2$.
The induced map $\gamma_j\colon \R/L\Z\to C_j$ is an
$m_j$-fold covering.
Hence both projections $\gamma_j\colon \R/L\Z\to S^2$ are embeddings
if and only if $m_1=m_2=1$, 
equivalently $\mu_1=\mu_2=2\pi/L$.
\end{proposition}

\begin{proof}
The vectors in \eqref{eq:circle-1}--\eqref{eq:circle-2} have unit norm.
Direct differentiation shows that their speeds are
$\cos\theta$ and $\sin\theta$, and that their signed curvature
coefficients are $\kappa_1$ and $\kappa_2$.  The congruence theorem
therefore proves the normal form.  Each component is a circle with
angular frequency $\mu_j$, so a common period exists exactly when the
ratio $\mu_1/\mu_2$ is rational.  
For the least positive period $L$ of $\gamma$, 
the integers $m_1,m_2$ in \eqref{eq:constant-data-period} are relatively prime,
and each induced map $\gamma_j\colon \R/L\Z\to C_j$ is an
$m_j$-fold covering.
\end{proof}

\subsection{K\"ahler magnetic geodesics}
\label{subsec:magnetic}

We conclude this section by applying the preceding curve theory to
K\"ahler magnetic geodesics.
We first write the magnetic equation for the product
K\"ahler--Einstein structure in terms of the two spherical projections.
We then classify its solutions by their curvature data and 
give an explicit normal form.

Let $J_{S^2}$ be the complex structure compatible with the standard
orientation and metric of $S^2$.  Thus
\[
 J_{S^2}X=x\times X
 \qquad
 (X\in T_xS^2).
\]
The product tensor fields
\begin{equation}\label{eq:positive-Kahler-structure}
 g=g_{S^2}\oplus g_{S^2},
 \qquad
 J=J_{S^2}\oplus J_{S^2}
\end{equation}
define a K\"ahler--Einstein structure on $M=S^2\times S^2$.
Its K\"ahler form is
\[
 \omega=g(J\,\cdot\,,\,\cdot\,)
 =\omega_{S^2}\oplus\omega_{S^2},
\]
where $\omega_{S^2}$ is the standard K\"ahler form on $S^2$.

A K\"ahler magnetic geodesic of charge $\lambda\in\R$ is a curve
$\gamma$ satisfying
\begin{equation}\label{eq:Kahler-magnetic-equation}
 \nabla_{\gamma'}\gamma'=\lambda J\gamma'.
\end{equation}
%Magnetic geodesics are regarded as nice perturbations of geodesics 
Magnetic geodesics generalize ordinary geodesics,
and have been studied very actively in dynamical systems, 
symplectic geometry as well as Riemannian geometry, 
see \cite{Arnold, Anosov-Sinai}. 
Equation~\eqref{eq:Kahler-magnetic-equation} implies that
every solution has constant speed.
In what follows, we restrict attention to unit-speed solutions.

\begin{theorem}\label{thm:Kahler-magnetic}
Every unit-speed K\"ahler magnetic geodesic of charge $\lambda$ in
$S^2\times S^2$ is a projection-regular curve with constant curvature data.
More precisely:
\begin{enumerate}
\item
in the G-projection-regular case,
$
 (\theta,\kappa_1,\kappa_2)
 =(\theta,\lambda,\lambda)
$
for some constant $\theta\in(0,\pi/2)$;
\item
in the E-projection-regular case, the signed curvature coefficient of the regular
projection is
$
 \kappa_0=\lambda.
$
\end{enumerate}
Conversely, every unit-speed projection-regular curve with
the corresponding curvature data in $(1)$ or $(2)$
is a K\"ahler magnetic geodesic of charge $\lambda$.
\end{theorem}

\begin{proof}
Equation~\eqref{eq:Kahler-magnetic-equation} splits into
\begin{equation}\label{eq:component-Kahler-magnetic}
 \nabla_s\gamma_j'
 =\lambda J_{S^2}\gamma_j',
 \qquad j=1,2.
\end{equation}
Taking the inner product with $\gamma_j'$ shows that each component
speed is constant.  Hence each component is either regular or constant,
so $\gamma$ is projection-regular.

Suppose first that $\gamma$ is G-projection-regular.  
Using \eqref{eq:FS}, we have
\[
 \nabla_s\gamma_1'
 =-\theta'\sin\theta\,\e_1
   +\kappa_1\cos\theta\,\n_1,\qquad
 \nabla_s\gamma_2'
 =\theta'\cos\theta\,\e_2
   +\kappa_2\sin\theta\,\n_2.
\]
Using $J_{S^2}\e_j=\n_j$,
comparison with \eqref{eq:component-Kahler-magnetic} gives
$
 \theta'=0$ and
$\kappa_1=\kappa_2=\lambda.
$
In the E-projection-regular case, the same calculation on the regular projection gives
$\kappa_0=\lambda$.  The converse follows immediately by substituting
these curvature data into the component equations.
\end{proof}

Combining Theorem~\ref{thm:Kahler-magnetic} with
Proposition~\ref{prop:projection-regular-homogeneous}, we obtain the following
consequence.

\begin{corollary}\label{cor:Kahler-magnetic-homogeneous}
Every K\"ahler magnetic geodesic in
$(S^2\times S^2,g,J)$ is homogeneous.
\end{corollary}

We next make the classification explicit.  For $c>0$, put
\begin{equation}\label{eq:Kahler-magnetic-component}
 \Gamma_{c,\lambda}(s)
 :=
 \frac1{\mu}
 \bigl(c\cos(\mu s),c\sin(\mu s),\lambda\bigr),
 \qquad
 \mu:=\sqrt{c^2+\lambda^2}.
\end{equation}
This is a spherical circle of speed $c$ and signed curvature
coefficient $\lambda$.

\begin{proposition}\label{prop:magnetic-normal-form}
Let $c_1,c_2$ be non-negative real numbers satisfying $c_1^2+c_2^2=1$.
For each $j\in \{1,2\}$, we set 
$$
\gamma_j(s):=\begin{cases}
\Gamma_{c_j,\lambda}(s) & c_j>0\\
X_j & c_j=0
\end{cases},
$$
where $X_j$ is an arbitrary fixed point in $S^2$.
Then the resulting curve
$
 \gamma(s)=\bigl(\gamma_1(s),\gamma_2(s)\bigr)
$
is a unit-speed K\"ahler magnetic geodesic of charge $\lambda$
in $(S^2\times S^2,g,J)$.
Every unit-speed K\"ahler magnetic geodesic is congruent to one of these
curves.
\end{proposition}

\begin{proof}
Direct differentiation of \eqref{eq:Kahler-magnetic-component} shows
that $\Gamma_{c,\lambda}$ has speed $c$ and satisfies
$
 \nabla_s\Gamma_{c,\lambda}'
 =\lambda J_{S^2}\Gamma_{c,\lambda}'.
$
The condition $c_1^2+c_2^2=1$ gives the unit-speed condition.
Conversely, Theorem~\ref{thm:Kahler-magnetic} shows that the component
speeds and curvature coefficients of any solution are constant.
The congruence theorem, applied to each nonconstant projection, then
gives the asserted normal form.
\end{proof}

\begin{remark}\label{rem:invariant-magnetic-extension}
K\"ahler magnetic geodesics for the neutral K\"ahler structure
\[
 g_{S^2}\oplus(-g_{S^2}),
 \qquad
 J_{S^2}\oplus(-J_{S^2})
\]
were studied by Godoy and Salvai \cite{Godoy-Salvai}.
This neutral K\"ahler surface is the mini-twistor space of $S^3$,
see, for example, \cite{BW}.
More generally, the componentwise argument used above applies to
arbitrary $K$-invariant metrics and invariant closed $2$-forms on
$S^2\times S^2$.
These tensor fields have the form
\[
 g_{p_1,p_2}
 =p_1g_{S^2}\oplus p_2g_{S^2},
 \qquad
 \omega_{q_1,q_2}
 =q_1\omega_{S^2}\oplus q_2\omega_{S^2},
\]
for some real numbers $p_1,p_2,q_1,q_2$
(see \cite{AGK,Anciaux}).
Our method used in Theorem~\ref{thm:Kahler-magnetic} 
yields a corresponding classification of their magnetic geodesics.
\end{remark}

\section{A Fenchel-type inequality for closed G-projection-regular curves in
$S^2\times S^2$}
\label{sec:fenchel-type}

This section studies a global invariant of closed G-projection-regular curves.
We first define the total curvature as the product Sasaki length of the
canonical lifts of the two spherical projections and compare it with
their ordinary total absolute curvatures.  We then prove a sharp
Fenchel-type inequality and determine its equality case.

\subsection{Total curvature}

In this subsection, we introduce the total curvature intrinsically from
the two spherical projections.  We first express it in terms of the
curvature triplet and then compare it with the Euclidean total absolute
curvatures of the projected curves.

Let
$
 \gamma=(\gamma_1,\gamma_2)\colon \R/L\mathbb Z
 \to S^2\times S^2
$
be a closed G-projection-regular curve, parametrized by arc length with respect to the
product metric $g$.  Put
\[
 v_1=\|\gamma_1'\|=\cos\theta,
 \qquad
 v_2=\|\gamma_2'\|=\sin\theta,
 \qquad
 \e_j=\frac{\gamma_j'}{v_j}.
\]
The canonical lift of the two projected curves is
\[
 \widehat\gamma
 :=\bigl((\gamma_1,\e_1),(\gamma_2,\e_2)\bigr)
 \colon \R/L\mathbb Z\longrightarrow US^2\times US^2.
\]

\begin{definition}\label{def:total-curvature-G-projection-regular}
The \emph{total curvature} of a closed G-projection-regular curve $\gamma$ is the length of
$\widehat\gamma$ with respect to the product Sasaki metric:
\begin{equation}\label{eq:total-curvature-G-projection-regular}
 \mathcal T(\gamma)
 :=\int_0^L\|\widehat\gamma'(s)\|_{\mathrm{Sas}}\,ds.
\end{equation}
\end{definition}

By the Frenet-type equations \eqref{eq:FS},
$
 \|\nabla_s\e_j\|=|\kappa_j|,
$
and therefore
\begin{equation}\label{eq:total-curvature-triplet}
 \mathcal T(\gamma)
 =\int_0^L
   \sqrt{1+\kappa_1(s)^2+\kappa_2(s)^2}\,ds.
\end{equation}
Although the right-hand side is written in the product arc-length
parameter, Definition~\ref{def:total-curvature-G-projection-regular} 
shows that the quantity is unchanged by a reparametrization.

\begin{proposition}\label{prop:total-curvature-comparison}
Let $\TC(\gamma_j)$ denote the Euclidean total absolute curvature of the
$j$th spherical projection, regarded as a space curve in $\R^3$.  Then
\begin{equation}\label{eq:total-curvature-comparison}
 \mathcal T(\gamma)\geq\sqrt{\TC(\gamma_1)^2+\TC(\gamma_2)^2}.
\end{equation}
Moreover, $\TC(\gamma_j)$ is the Sasaki length of the canonical lift
$(\gamma_j,\e_j)$.
\end{proposition}

\begin{proof}
The speed of $(\gamma_j,\e_j)$ in the Sasaki metric is
\[
 A_j:=\sqrt{\|\gamma_j'\|^2+\|\nabla_s\e_j\|^2}
 =\sqrt{v_j^2+\kappa_j^2}.
\]
This is also the Euclidean curvature density of the spherical curve,
because its normal curvature in the unit sphere is one.  
The inequality in \eqref{eq:total-curvature-comparison} 
is the integral triangle inequality for vectors.
\end{proof}

\subsection{The inequality and its equality case}

We now combine the preceding comparison with the classical Fenchel
theorem.  We first derive the lower bound, then analyze equality in both
the integral triangle inequality and Fenchel's theorem, and finally
record the values attained by closed constant-data examples.

\begin{theorem}\label{prop:fenchel-triplet}
Let $\gamma\colon \R/L\mathbb Z\to S^2\times S^2$ be a closed unit-speed
G-projection-regular curve.  Then
\begin{equation}\label{eq:fenchel-triplet}
 \mathcal T(\gamma)\geq2\sqrt2\,\pi.
\end{equation}
Equality holds if and only if both 
projections $\gamma_j\colon \R/L\Z\to S^2$ $(j=1,2)$ are embeddings
and the curvature triplet $(\theta,\kappa_1,\kappa_2)$ is constant.
\end{theorem}

\begin{proof}
Let $u_j$ be the arc-length parameter of $\gamma_j$ and let $k_{g,j}$ be
its signed geodesic curvature in $S^2$.  Since $du_j=v_jds$, equation
\eqref{eq:FS} gives
\[
 \frac{d\e_j}{du_j}
 =-\gamma_j+\frac{\kappa_j}{v_j}\n_j,
 \qquad
 k_{g,j}=\frac{\kappa_j}{v_j}.
\]
Regarded as a space curve in $\R^3$, $\gamma_j$ has Euclidean curvature
$k_{E,j}=\sqrt{1+k_{g,j}^2}$.  Hence its total absolute curvature is
\begin{equation}\label{eq:component-total-curvature}
 \TC(\gamma_j)
 =\int_0^L A_j(s)\,ds,
 \qquad
 A_j:=\sqrt{v_j^2+\kappa_j^2}.
\end{equation}
Since
$
 A_1^2+A_2^2=1+\kappa_1^2+\kappa_2^2,
$
the triangle inequality for the $\R^2$-valued function $(A_1,A_2)$ gives
\begin{align}
\mathcal T(\gamma)
 =\int_0^L\sqrt{A_1^2+A_2^2}\,ds
 &\geq
 \sqrt{\left(\int_0^LA_1\,ds\right)^2
       +\left(\int_0^LA_2\,ds\right)^2}\notag\\
 &=\sqrt{\TC(\gamma_1)^2+\TC(\gamma_2)^2}.
\label{eq:Minkowski-Fenchel}
\end{align}
Fenchel's theorem yields
$\TC(\gamma_j)\geq2\pi$ for $j=1,2$, proving
\eqref{eq:fenchel-triplet}.

Suppose that equality holds.  Equality is then required in both steps, so
\begin{equation}\label{eq:component-Fenchel-equality}
 \TC(\gamma_1)=\TC(\gamma_2)=2\pi.
\end{equation}
The equality case of Fenchel's theorem implies that 
each $\gamma_j$ is an embedding whose image is a convex planar curve.
Since the image lies on the unit sphere, 
it is a plane section of $S^2$, 
hence a spherical circle. 
Its signed geodesic curvature is therefore constant.

Equality in \eqref{eq:Minkowski-Fenchel} implies that the vector
$(A_1(s),A_2(s))$
has a fixed direction.
Thus $A_1/A_2$ is constant.
Together with \eqref{eq:component-Fenchel-equality}, this gives
$A_1=A_2$.  Since
$
 A_j=v_j\sqrt{1+k_{g,j}^2}
$
and both $k_{g,j}$ are constant, the ratio $v_1/v_2$ is constant.  The
relations $v_1^2+v_2^2=1$ and $v_1,v_2>0$ now imply that both component
speeds, and hence $\theta$, are constant.  Finally,
$\kappa_j=v_jk_{g,j}$ is constant.

Conversely, assume that both 
$\gamma_1$ and $\gamma_2$ are embeddings
and that the curvature triplet is constant.
Then the image of each $\gamma_j$ is a spherical circle,
and $A_j$ is constant.
Therefore
\[
 LA_j=\TC(\gamma_j)=2\pi
 \qquad (j=1,2).
\]
Consequently,
$
 \mathcal T(\gamma)
 =L\sqrt{A_1^2+A_2^2}
 =2\sqrt2\,\pi,
$
as required.
\end{proof}

The assumption that both $\gamma_1$ and $\gamma_2$ are embeddings 
is essential in the converse.
For a closed G-projection-regular curve with constant curvature data,
each $\gamma_j$ is an $m_j$-fold covering of its image circle
for some positive integer $m_j$, and $\TC(\gamma_j)=2\pi m_j$.
Equality in \eqref{eq:fenchel-triplet} excludes $m_j>1$, 
as illustrated by the following example.

\begin{example}[Closed G-projection-regular curves with constant curvature triplet]
\label{ex:total-curvature-constant-data}
Let $\gamma\colon \R/L\Z\to S^2\times S^2$ be a closed G-projection-regular curve with
constant curvature triplet
$(\theta,\kappa_1,\kappa_2)$.
Put
\[
 \mu_1=\sqrt{\cos^2\theta+\kappa_1^2},
 \qquad
 \mu_2=\sqrt{\sin^2\theta+\kappa_2^2}.
\]
Since the projected curves are spherical circles with angular
frequencies $\mu_1$ and $\mu_2$, their closedness implies that
there exist positive integers $m_1,m_2$ such that
\[
 \mu_1L=2\pi m_1,
 \qquad
 \mu_2L=2\pi m_2.
\]
Here, each $\gamma_j$ is an $m_j$-fold covering
of its image circle $C_j=\gamma_j(\R)$.

Since
$
 1+\kappa_1^2+\kappa_2^2
 =\mu_1^2+\mu_2^2,
$
we obtain
$
 \mathcal T(\gamma)
 =2\pi\sqrt{m_1^2+m_2^2}.
$
Thus
$
 \mathcal T(\gamma)\ge 2\sqrt{2}\,\pi,
$
and equality holds if and only if
$
 m_1=m_2=1,
$
or equivalently, if and only if both 
$\gamma_1$ and $\gamma_2$ are embeddings.
\end{example}

\section{From curves in $S^2\times S^2$ to ruled surfaces in $S^3$}
\label{sec:LS3}

This section transfers the curve theory developed above to ruled
surfaces in $S^3$.  We first identify
$S^2\times S^2$ equivariantly with the space of oriented geodesics in
$S^3$ and give an explicit formula for those geodesics.  
We next introduce ruling Gauss maps of a ruled surface.
We express the singularities and curvatures of ruled surfaces 
through the curvature data of its associated curve. 
We then characterize helicoidal ruling families and 
describe those arising from K\"ahler magnetic geodesics.
Finally, we apply the Fenchel-type inequality to closed G-projection-regular ruled
surfaces and analyze its equality case.

\subsection{The space of oriented geodesics in $S^3$}
\label{subsec:oriented-geodesic-space}

In this subsection, we first construct the product realization of the
space $\LS$ of oriented geodesics by means of the Hodge decomposition.
We then compare the induced isometry action with the factorwise action
used in the preceding sections.  Finally, we introduce affine
coordinates and write the corresponding oriented geodesics explicitly.

An oriented geodesic in $S^3$ is the intersection of $S^3$ with an
oriented $2$-plane in $\R^4$.  Hence
\[
 \LS=\widetilde{\operatorname{Gr}}_2(\R^4).
\]
If $(u,v)$ is a positively oriented orthonormal basis of such a plane,
put $p=u\wedge v$.  With respect to the Hodge decomposition
\[
 \Lambda^2\R^4=\Lambda^2_+\R^4\oplus\Lambda^2_-\R^4,
 \qquad
 p_\pm=\frac12(p\pm *p),
\]
one has $\|p_+\|=\|p_-\|=1/\sqrt2$.  The oriented plane therefore
determines
\[
 \bigl(\sqrt2p_+,\sqrt2p_-\bigr)
 \in S(\Lambda^2_+\R^4)\times S(\Lambda^2_-\R^4),
\]
where $S(\Lambda^2_{\pm}\R^4)$ denotes the unit sphere in
$\Lambda^2_{\pm}\R^4$.
Conversely, if $(x_+,x_-)$ belongs to this product, then
$p=(x_++x_-)/\sqrt2$ has norm one and
\[
 p\wedge p
 =\frac12\bigl(\|x_+\|^2-\|x_-\|^2\bigr)
   \operatorname{vol}_{\R^4}=0.
\]
Thus $p$ is decomposable and represents a unique oriented $2$-plane.
After choosing orientation-preserving Euclidean identifications
$\Lambda_\pm^2\R^4\simeq\R^3$, we obtain
\begin{equation}\label{eq:LS-product}
 \LS\cong\widetilde{\operatorname{Gr}}_2(\R^4)
 \cong S^2\times S^2.
\end{equation}

The induced homomorphism
\[
 \SO(4)\longrightarrow
 \SO(\Lambda^2_+\R^4)\times\SO(\Lambda^2_-\R^4)
 \cong\SO(3)\times\SO(3)=K
\]
is onto and has kernel $\{\pm\operatorname{Id}\}$.
Consequently, \eqref{eq:LS-product} is equivariant and 
the isometric action on $\LS$ induced by $\SO(4)\cong\Isom_0(S^3)$ 
is exactly the factorwise action used in
Sections~\ref{sec:congruence-thm} and \ref{sec:fenchel-type}.

For explicit formulas, we use the matrix model
\[
 \H=\left\{
 \begin{pmatrix}
 x_0+ix_3&-x_2+ix_1\\
 x_2+ix_1&x_0-ix_3
 \end{pmatrix}:
 \ x_0,x_1,x_2,x_3\in\R
 \right\}
\]
of the skew field of quaternions.
Then $S^3$ is identified with $\SU(2)$, and the unit sphere in
$\Im\H\simeq\R^3$ is identified with $S^2$.  The spin covering
\[
 \SU(2)\times\SU(2)\longrightarrow \SO(4)
\]
is induced by $(A,B)\cdot X=AXB^{-1}$.

Via the quadric model, the identification \eqref{eq:LS-product} may be
written as
\[
 \LS\cong Q^2\cong\widehat\C\times\widehat\C.
\]
Let $\mathcal U\subset\LS$ be the open dense subset corresponding to the
standard affine chart
$
 \C^2=(\widehat\C\setminus\{\infty\})
       \times(\widehat\C\setminus\{\infty\}).
$
The following formula fixes the normalization of these coordinates.

\begin{proposition}\label{prop:geodesic-representation}
If $\ell\in\mathcal U$ has affine coordinates $(z_1,z_2)\in\C^2$, then
$\ell=[G^{z_1,z_2}]$, where
\begin{equation}\label{eq:Gamma-z1z2}
 G^{z_1,z_2}(t)
 =P(z_1)L(t)P(\overline{z_2})^{-1},
 \qquad t\in\R,
\end{equation}
and
\begin{equation}\label{eq:matrix}
 P(z)=\frac1{\sqrt{1+|z|^2}}
 \begin{pmatrix}1&z\\-\bar z&1\end{pmatrix},
 \qquad
 L(t)=\begin{pmatrix}e^{it}&0\\0&e^{-it}\end{pmatrix}.
\end{equation}
\end{proposition}

\begin{proof}
For background on the product model
$Q^2\cong\widehat\C\times\widehat\C$, see \cite{Fujimoto}.
With the affine-coordinate normalization used here, the corresponding
geodesic is represented by
\begin{equation}\label{eq:Gamma-expanded}
\frac1{\sqrt{(1+|z_1|^2)(1+|z_2|^2)}}
\begin{pmatrix}
 e^{it}+e^{-it}z_1z_2&
 e^{-it}z_1-e^{it}\bar z_2\\
 -e^{it}\bar z_1+e^{-it}z_2&
 e^{-it}+e^{it}\bar z_1\bar z_2
\end{pmatrix}.
\end{equation}
Both $P(z_1)$ and $P(\overline{z_2})$ belong to $\SU(2)$,
so \eqref{eq:Gamma-z1z2} is an oriented unit-speed geodesic, being
an isometric image of $L(t)$.
A direct multiplication gives \eqref{eq:Gamma-expanded}.
\end{proof}

Define $\Psi\colon \C^2\to\mathcal U$ by
\begin{equation}\label{eq:Psi}
 \Psi(z_1,z_2)=[G^{z_1,z_2}].
\end{equation}
For an oriented unit-speed geodesic $G$, put
\[
 \psi_R([G])=-\dot GG^{-1},\qquad
 \psi_L([G])=-G^{-1}\dot G.
\]
Both expressions are constant along $G$, because the metric on
$S^3=\SU(2)$ is bi-invariant, and they are unchanged by a parameter
translation.  Hence $\psi_R,\psi_L\colon \LS\to S^2$ are well defined.
Let $\Pi\colon S^2\to S^2$ be the reflection
\begin{equation}\label{eq:reflection}
 \Pi
 \begin{pmatrix}
  ix_3&-x_2+ix_1\\
  x_2+ix_1&-ix_3
 \end{pmatrix}
 =
 \begin{pmatrix}
  ix_3&x_2+ix_1\\
  -x_2+ix_1&-ix_3
 \end{pmatrix}.
\end{equation}
Note that $\Pi$ is just a transpose of a matrix.
Let $\pi_N$ be stereographic projection from the north pole:
\begin{equation}\label{eq:stereoproj-N}
 \pi_N
 \begin{pmatrix}
  ix_3&-x_2+ix_1\\
  x_2+ix_1&-ix_3
 \end{pmatrix}
 =\frac{x_1+ix_2}{1-x_3}.
\end{equation}

\begin{proposition}\label{prop:psi-inverse}
The inverse of $\Psi$ is
\begin{equation}\label{eq:psi-inverse}
 \Psi^{-1}
 =(\pi_N\circ\psi_R,\,\pi_N\circ\Pi\circ\psi_L)\colon
 \mathcal U\longrightarrow\C^2.
\end{equation}
\end{proposition}

\begin{proof}
Substitution of \eqref{eq:Gamma-z1z2} gives
\[
 \pi_N\circ\psi_R([G^{z_1,z_2}])=z_1,\qquad
 \pi_N\circ\Pi\circ\psi_L([G^{z_1,z_2}])=z_2.
\]
Thus the right-hand side of \eqref{eq:psi-inverse} is a left inverse of
$\Psi$.  Proposition~\ref{prop:geodesic-representation} says that
$\Psi$ is onto $\mathcal U$, so it is the inverse.
\end{proof}

\subsection{Ruled surfaces and ruling Gauss maps}
\label{subsec:ruled-gauss}

We now pass from individual oriented geodesics to one-parameter
families.  We first define ruled surfaces and show that, up to
congruence, their associated curves lie in the affine chart introduced
above.  We then define the right and left ruling Gauss maps and identify
them with the two spherical components of the associated curve.

\begin{definition}\label{def:ruled-surface}
A \emph{ruled surface} in $S^3$ is a smooth one-parameter family
$\{\ell_s\}_{s\in I}$ of oriented geodesics such that
the corresponding curve $\gamma\colon I\to\LS$, defined by
$\gamma(s)=\ell_s$, is regular.
Each $\ell_s$ is called an \emph{oriented ruling} of the surface,
and $\gamma$ is called its \emph{associated curve}.
\end{definition}

Let $G_s(t)$ be a unit-speed parametrization of $\ell_s=[G_s]$.
Then 
$$
F \colon I \times S^1 \to S^3,\qquad
F(s,t):=G_s(t)
$$ 
gives a parametrization of a ruled surface.  
Note that $F$ is not assumed to be an immersion. 

\begin{proposition}[Representation formula for ruled surfaces]%
\label{prop:ruled-surface}
Up to congruence, every ruled surface admits a parametrization
\begin{equation}\label{eq:ruled-representation}
 F(s,t)=P(z_1(s))L(t)P(\overline{z_2(s)})^{-1}
 \qquad (s\in I,\ t\in\R/2\pi\mathbb Z),
\end{equation}
where $z=(z_1,z_2)\colon I\to\C^2$ is regular.
\end{proposition}

\begin{proof}
Each projected curve omits a point of $S^2$ by Sard's theorem.  The
factorwise $\SO(3)\times\SO(3)$ action moves these two points to the
north poles.  Both projections then lie in the standard affine charts,
and \eqref{eq:Gamma-z1z2} gives the formula.
\end{proof}

\begin{definition}\label{def:ruling-Gaussmap}
For a ruled surface $F\colon I\times S^1\to S^3=\SU(2)$, 
define
\begin{equation}\label{eq:ruling-Gauss-maps}
 \nu_R(s)=-F_t(s,t)F(s,t)^{-1},\qquad
 \nu_L(s)=-F(s,t)^{-1}F_t(s,t).
\end{equation}
The geodesic equation and the bi-invariance of the metric show that
these expressions are independent of $t$.  They are the \emph{right}
and \emph{left ruling Gauss maps}.  
\end{definition}

The reflection $\Pi$ in
\eqref{eq:reflection} relates them to the product coordinates as follows.

\begin{proposition}\label{prop:ruling-Gauss-curve}
The associated curve $\gamma$ of a ruled surface is given by
\begin{equation}\label{eq:Gauss-curve-identification}
 \gamma=(\nu_R,\Pi\circ\nu_L)\colon I\longrightarrow S^2\times S^2.
\end{equation}
\end{proposition}

\begin{proof}
For every fixed $s$, the two expressions in
\eqref{eq:ruling-Gauss-maps} are precisely
$\psi_R([G_s])$ and $\psi_L([G_s])$.  The conclusion is therefore
Proposition~\ref{prop:psi-inverse}, first on $\mathcal U$ and then on all
of $\LS$ by continuity.
\end{proof}

\subsection{Projection-regular ruled surfaces and curvature}
\label{subsec:projection-regular-ruled}

This subsection translates the local curvature data of a projection-regular curve into
the geometry of the corresponding ruled surface.  
We first define projection-regular,
G-projection-regular, and E-projection-regular ruled surfaces 
and determine their singular sets.

\begin{definition}\label{def:projection-regular-ruled-surface}
Let $\gamma$ be the associated curve of a ruled surface.
We call the ruled surface {\it projection-regular} 
if $\gamma$ is projection-regular.
It is {\it G-projection-regular} or {\it E-projection-regular} 
according as $\gamma$ is generic or exceptional.
We define the curvature data of the ruled surface 
to be the curvature data of its associated curve.
\end{definition}

\begin{proposition}\label{prop:singular-set-projection-regular}
Let $F$ be a projection-regular ruled surface with unit-speed associated curve.
\begin{enumerate}
\item An E-projection-regular ruled surface has no singular points.
\item For a G-projection-regular ruled surface, one can choose the origins of the ruling
parameters locally so that
\begin{equation}\label{eq:det-first-form}
 \det(g_{ij})(s,t)=\frac14W(s,t),\qquad
 W(s,t):=1-\sin2\theta(s)\cos2t,
\end{equation}
where $g_{ij}$ are the metric coefficients of $F$.
Thus a ruling contains singular points if and only if
$\theta(s)=\pi/4$; in the adapted parameter they occur at
$t=0,\pi$ modulo $2\pi$.
\end{enumerate}
\end{proposition}

\begin{proof}
At a fixed $s$, use a congruence to put both affine coordinates at the
origin and both tangent directions on the positive real axes.
Put
\[
 a=\frac12\|\gamma_1'\|,\qquad
 b=\frac12\|\gamma_2'\|.
\]
At the fixed point, differentiation of
\eqref{eq:ruled-representation} and right translation to the identity
give
\[
 F_sF^{-1}
 =(-a+b\cos2t)
   \begin{pmatrix}0&-1\\1&0\end{pmatrix}
 -b\sin2t
   \begin{pmatrix}0&i\\i&0\end{pmatrix},
\quad
 F_tF^{-1}
 =\begin{pmatrix}i&0\\0&-i\end{pmatrix}.
\]
These three matrices are orthonormal in $\Im\H$.  Therefore
\begin{equation}\label{eq:first-form-adapted}
 g_{11}=a^2+b^2-2ab\cos2t,\qquad
 g_{12}=0,\qquad
 g_{22}=1.
\end{equation}
Equivalently,
$
 g_{11}g_{22}-g_{12}^2
 =(1-2v_1v_2\cos2t)/4,
$
where $v_j=\|\gamma_j'\|$.  For a G-projection-regular curve,
$(v_1,v_2)=(\cos\theta,\sin\theta)$,
and for an E-projection-regular curve the speeds are
$(1,0)$ or $(0,1)$.  The assertions follow.
\end{proof}

Surfaces of constant Gaussian curvature or constant mean curvature in
$\S^3$ form important classes, and some of them are realized as ruled
surfaces.  We first consider ruled surfaces of constant mean curvature.
It is known that, up to congruence and locally on the regular set,
every non-totally-geodesic minimal ruled surface in $\S^3$
is parametrized, for some $h\in\R$, by
\begin{equation}\label{eq:Lawson-helicoid}
F(s,t)
=
(\cos s\cos t,\,
 \cos(hs)\sin t,\,
 \sin(hs)\sin t,\,
 \sin s\cos t).
\end{equation}
Such a surface is called a \emph{Lawson helicoid}
\cite[Proposition~7.2]{Lawson}.
Related minimal surfaces in $S^3$ foliated by circles were studied by
Kutev and Milousheva \cite{Kutev-Milousheva}.

Moreover, locally on the regular set, every nonminimal ruled surface
of constant mean curvature in $\S^3$ is, up to congruence,
an open subset of a product torus
\[
S^1(r)\times S^1\bigl(\sqrt{1-r^2}\bigr)\subset S^3,
\qquad
0<r<1,\quad r\ne\frac1{\sqrt2},
\]
as follows from \cite[Theorem~3]{Park2019}; see also
\cite[Theorem~1.4]{Barbosa-Delgado} for a classification under a
completeness assumption. These classification results yield the following.

\begin{proposition}\label{prop:CMC}
The associated curve of
a non-totally-geodesic minimal ruled surface in $S^3$
is a geodesic in $\LS$, when parametrized by arc length.
In the G-projection-regular case, its curvature data satisfy
$\theta'=0$ and $\kappa_1=\kappa_2=0$, and
in the E-projection-regular case, $\kappa_0=0$.
Moreover, the curve in $\LS$ corresponding to a nonminimal CMC ruled
surface in $S^3$ is an E-projection-regular curve of the form
\[
\gamma(s)=(X,\gamma_0(s))
\quad\text{or}\quad
\gamma(s)=(\gamma_0(s),X),
\]
where $X\in S^2$ is constant and $\gamma_0$ is a small circle in $S^2$.
In particular, all these curves are projection-regular.
\end{proposition}

For totally geodesic minimal ruled surfaces, see
Proposition~\ref{prop:CMC-CGC}.

\begin{proof}
The Lawson helicoid \eqref{eq:Lawson-helicoid} can be written in
$\SU(2)$ as
\[
F(s,t)
=
\begin{pmatrix}
e^{is}\cos t & i e^{ihs}\sin t\\
i e^{-ihs}\sin t & e^{-is}\cos t
\end{pmatrix}.
\]
By Proposition~\ref{prop:psi-inverse}, 
the corresponding curve
$\gamma=(\gamma_1,\gamma_2)$ is given by
$\gamma_1=-F_tF^{-1}$ and $\gamma_2=-\Pi(F^{-1}F_t)$,
where $\Pi\colon S^2\to S^2$ is the reflection defined in
\eqref{eq:reflection}.
By direct calculation, we have
\[
\gamma_1(s)
=\begin{pmatrix}
 0 & -i e^{i (h+1) s} \\
 -i e^{-i (h+1) s} & 0 \\
\end{pmatrix},\quad
\gamma_2(s)
=\begin{pmatrix}
 0 & -i e^{-i (h-1) s} \\
 -i e^{i (h-1) s} & 0 \\
\end{pmatrix}.
\]
Hence, $\gamma$ is a geodesic in $\LS$.

Similarly, for $0<\beta<\pi/2$ with $\beta\ne\pi/4$,
a nonminimal parallel torus of the Clifford torus can be written in
$\SU(2)$ as
\[
F(s,t)
=
\begin{pmatrix}
\cos\beta\,e^{it}&-\sin\beta\,e^{i(t+s)}\\
\sin\beta\,e^{-i(t+s)}&\cos\beta\,e^{-it}
\end{pmatrix}.
\]
The assertion now follows from a calculation analogous to the one
above.
\end{proof}

\begin{proposition}\label{prop:CMC-CGC}
Let $\gamma\colon  I\to S^2\times S^2$ be a regular curve,
and  
$F \colon I\times \S^1 \to \S^3$ be the ruled surface given by $\gamma$.
If the Gaussian curvature $K$ of $F$ is constant
on the regular set, then 
either $K=0$ $($flat$)$ or $K=1$ $($extrinsically flat$)$ holds.
Moreover the following assertions hold.
\begin{enumerate}
\item
If $F$ is flat $(K=0)$,
then $\gamma$ is an exceptional projection-regular curve.
\item
If $F$ is extrinsically flat $(K=1)$,
then $\gamma$ is a generic projection-regular curve
whose angle function satisfies $\theta=\pi/4$.
\item
If $F$ is totally geodesic,
$\gamma$ is a generic projection-regular curve
whose angle function satisfies $\theta=\pi/4$
and $\kappa_1=\kappa_2$.
\end{enumerate}
\end{proposition}

\begin{proof}
  Let $z=(z_1,z_2)\colon I\to \C^2$ be a regular curve.
The ruled surface given by $z(s)$
is parametrized as
$$
F(s,t) = P(z_1(s)) L(t) P(\overline{z_2(s)})^{-1}
\qquad
(s\in I,\, t\in \S^1)
$$ 
as in Proposition~\ref{prop:ruled-surface},
where $P(z)$, $L(t)$ are the matrices given in \eqref{eq:matrix}.
We fix a point $s_0\in I$ arbitrarily.
Using an isometry of $\S^3$,
we may assume that 
\[
  z_1(s_0)=z_2(s_0)=0,\qquad
  z_1'(s_0)=a,\qquad 
  z_2'(s_0)=b,
\]
where $a$, $b$ are non-negative real numbers satisfying $(a,b)\ne(0,0)$.
Then, the Gaussian curvature $K$ and the mean curvature $H$ 
at $(s_0,t)$ are respectively written as
\begin{align}
\label{eq:Gaussian-F}
  K(s_0,t)
  &= 1-\frac{(a^2-b^2)^2}{\left(a^2+b^2 - 2 a b \cos 2 t\right)^2},\\
\label{eq:mean-F}
  H(s_0,t)
  &=\frac{
         -a\Im \left(e^{-2i t} z_2''(s_0)+ z_1''(s_0)\right)
         + b\Im \left(e^{-2i t} z_1''(s_0)+ z_2''(s_0)\right)
  }{2\left(\sqrt{a^2+b^2 - 2 ab \cos 2t }\right)^3}.
\end{align}
Suppose that $K$ is constant.  For a fixed $s_0$,
\eqref{eq:Gaussian-F} is independent of $t$ only if $a=b$ or $ab=0$.
Since the value of $K$ is independent of $s_0$, there are precisely two
possibilities.  If $a=b$ at every point, then $K=1$;
hence $z$ is a 
G-projection-regular curve
with $\theta\equiv\pi/4$.  If $ab=0$
at every point, then $K=0$;
hence $z$ is an exceptional projection-regular curve. Conversely, \eqref{eq:Gaussian-F}
shows that these two cases give $K=1$ and $K=0$, respectively.  This
proves (1) and (2).

We next prove (3).  If $F$ is totally geodesic, then $K=1$ and
$H=0$.  Thus $a=b$, and \eqref{eq:mean-F} gives
$z_1''(s_0)=z_2''(s_0)$.
Therefore, 
$\kappa_1(s_0)=\kappa_2(s_0)$, proving (3).
\end{proof}

\begin{remark}
Ruled flat surfaces are called Hopf cylinders \cite{Pinkall}.
General closed flat surfaces, not necessarily ruled, 
were classified by Kitagawa \cite{Kitagawa}.
Although every closed immersed extrinsically flat surface is totally
geodesic \cite{ONeil-Stiel}, numerous nontrivial examples exist
when singularities are allowed.
The classification of closed extrinsically flat surfaces 
with wave front singularities is given in \cite{Honda2017}.
\end{remark}

\subsection{Helicoidal and K\"ahler magnetic ruling families}
\label{subsec:helicoidal}

We next interpret constant curvature data through ambient symmetries of
$S^3$.  We first characterize helicoidal ruling families as those whose
associated curves have constant curvature data.  We then describe the
closed constant-data examples and apply the classification of
K\"ahler magnetic geodesics from
Subsection~\ref{subsec:magnetic}.

\begin{definition}\label{def:helicoidal-motion}
A one-parameter subgroup $\{\Phi_s\}$ of
$\Isom_0(S^3)\cong\SO(4)$ is called
a \emph{helicoidal motion}.
A ruled surface with specified oriented rulings
$\{\ell_s\}$ is \emph{helicoidal} if
\[
 \ell_s=\Phi_s(\ell_0).
\]
\end{definition}

After an orientation-preserving congruence, every helicoidal motion has
the following double-rotation form:
\begin{equation}\label{eq:double-rotation}
 \Phi_s(\zeta,\eta)
 =\bigl(e^{ias}\zeta,e^{ibs}\eta\bigr),
 \qquad a,b\in\R,
\end{equation}
on $S^3\subset\C^2$.
We say that a nontrivial helicoidal motion 
is a \emph{rotational motion} if all its elements
fix the same great circle in $S^3$ pointwise.
For the motion in \eqref{eq:double-rotation}, this is equivalent to
$ab=0$ and $(a,b)\ne(0,0)$.
The angular speeds of the two nonconstant spherical projections are
\begin{equation}\label{eq:projection-frequencies}
 \mu_1=|a+b|,\qquad \mu_2=|a-b|.
\end{equation}

\begin{proposition}\label{prop:helicoidal-characterization}
For a ruled surface with unit-speed associated curve $\gamma\colon \R\to\LS$,
the following are equivalent:
\begin{enumerate}
\item $\gamma$ is projection-regular with constant curvature data;
\item $\gamma$ is homogeneous;
\item the corresponding ruled surface is helicoidal.
\end{enumerate}
\end{proposition}

\begin{proof}
The equivalence of (1) and (2) is
Proposition~\ref{prop:projection-regular-homogeneous}.  If
$\gamma(s)=\Phi_s\cdot\gamma(0)$, equivariance of the correspondence
between points of $\LS$ and oriented geodesics gives
$\ell_s=\Phi_s(\ell_0)$, and the converse is immediate.
\end{proof}

\begin{example}[G-projection-regular ruled surfaces with constant curvature data]
\label{ex:G-projection-regular-ruled-constant-data}
Let
$
 \gamma\colon \R/L\Z\to\LS
$
be a closed G-projection-regular curve with constant curvature triplet
$(\theta,\kappa_1,\kappa_2)$, as considered in
Example~\ref{ex:total-curvature-constant-data}.
Let
$
 F\colon T_L^2\to S^3
$
be the ruled surface corresponding to $\gamma$,
where 
$T_L^2 = (\mathbb{R}/L \mathbb{Z}) \times (\mathbb{R}/2\pi \mathbb{Z})$.
Then $F$ is invariant under a one-parameter subgroup of
$\operatorname{Isom}_0(S^3)$ and hence is a helicoidal ruled
surface.
\end{example}

\begin{corollary}\label{cor:magnetic-ruled-invariant}
Every ruled surface whose associated curve corresponds, under
\eqref{eq:LS-product}, to a K\"ahler magnetic geodesic of charge
$\lambda$ is helicoidal.  Its curvature data are
\[
 (\theta,\kappa_1,\kappa_2)=(\theta,\lambda,\lambda)
 \quad\text{or}\quad
 \kappa_0=\lambda.
\]
\end{corollary}

\begin{figure}[htb]
\begin{center}
 \begin{tabular}{c@{\hspace{8mm}}c}
  \resizebox{4.5cm}{!}{\includegraphics{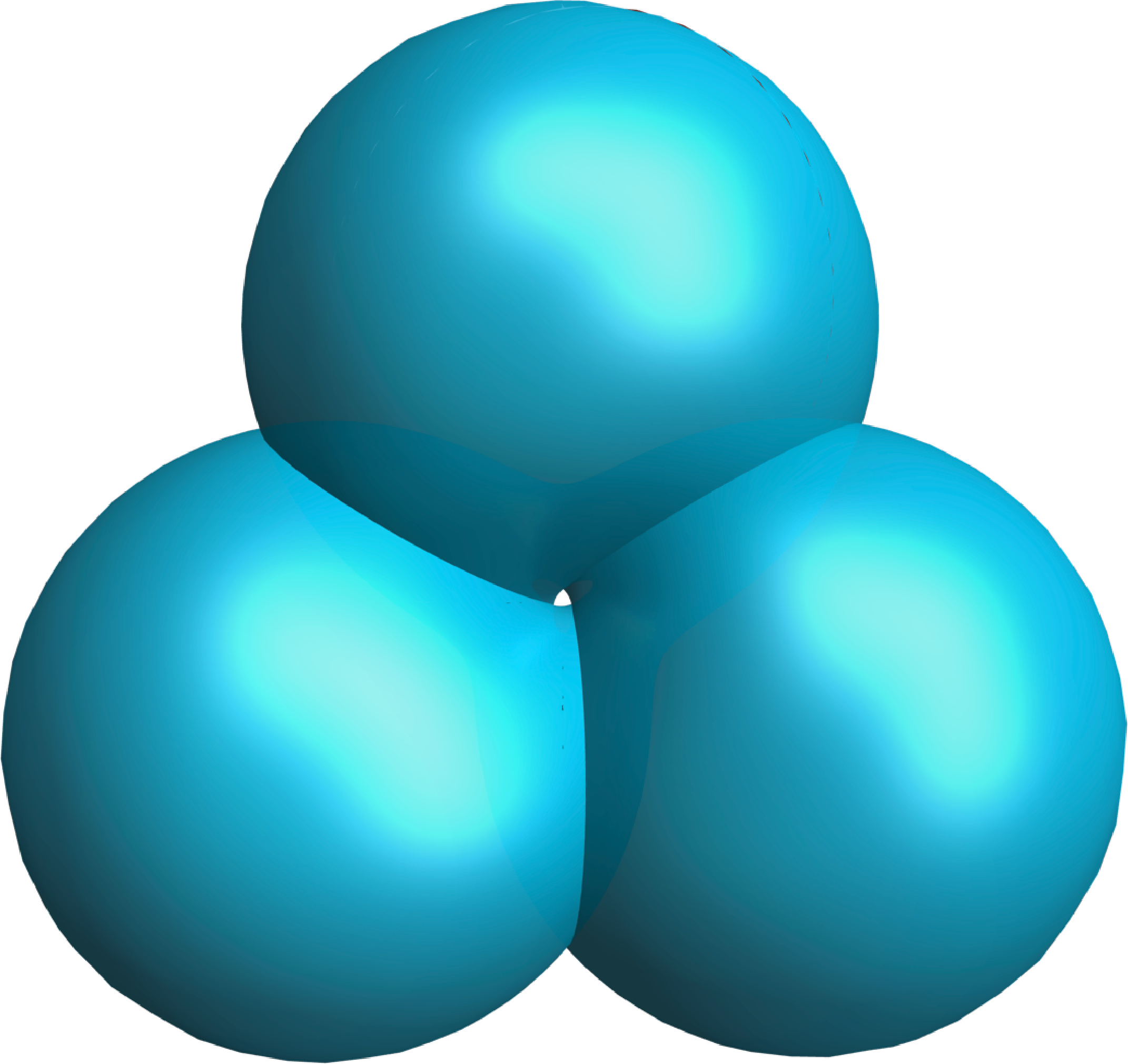}}&
  \resizebox{4.5cm}{!}{\includegraphics{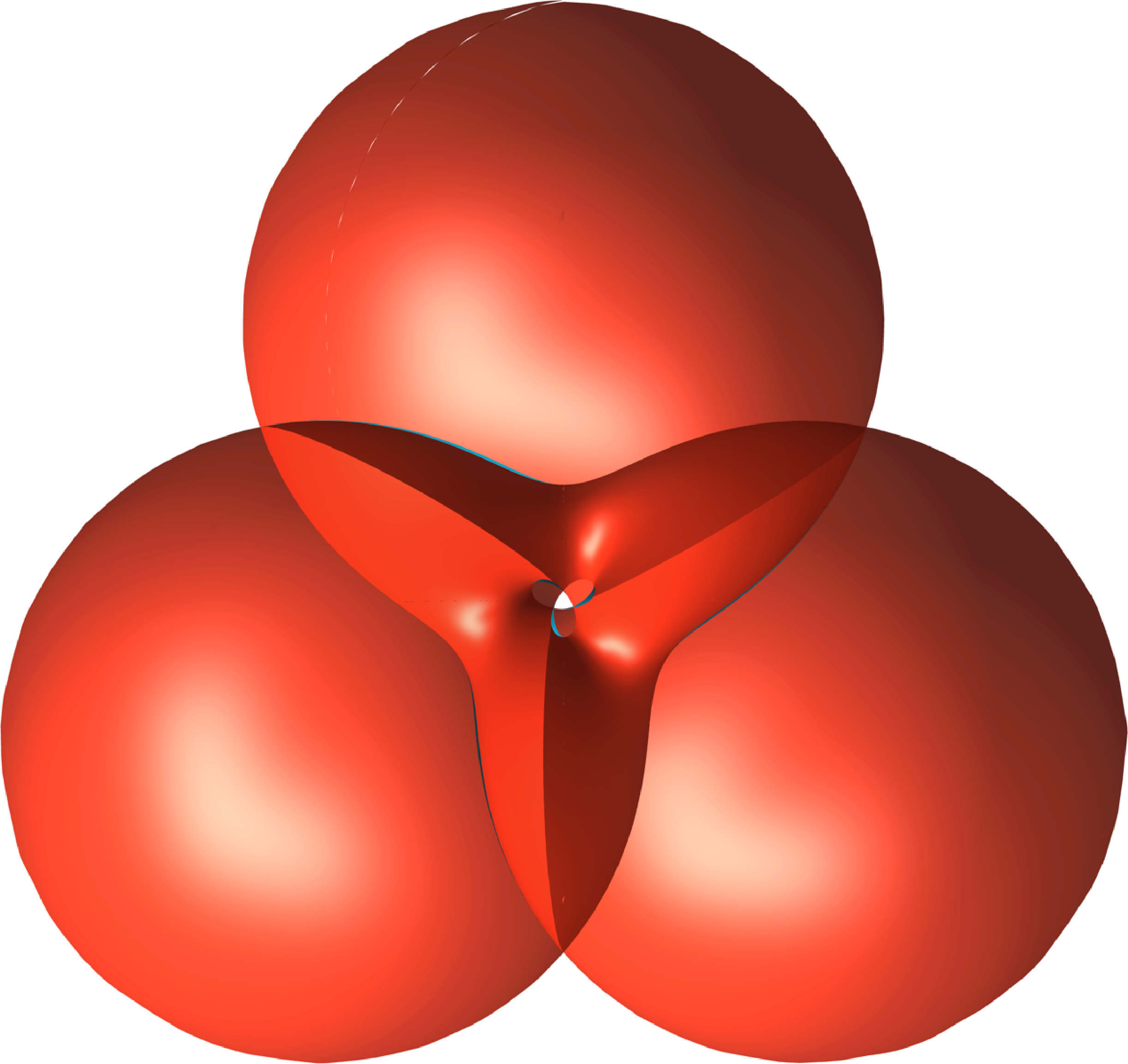}}\\
  {\footnotesize (a1)} & 
  {\footnotesize (a2)} \\ 
  \resizebox{4.5cm}{!}{\includegraphics{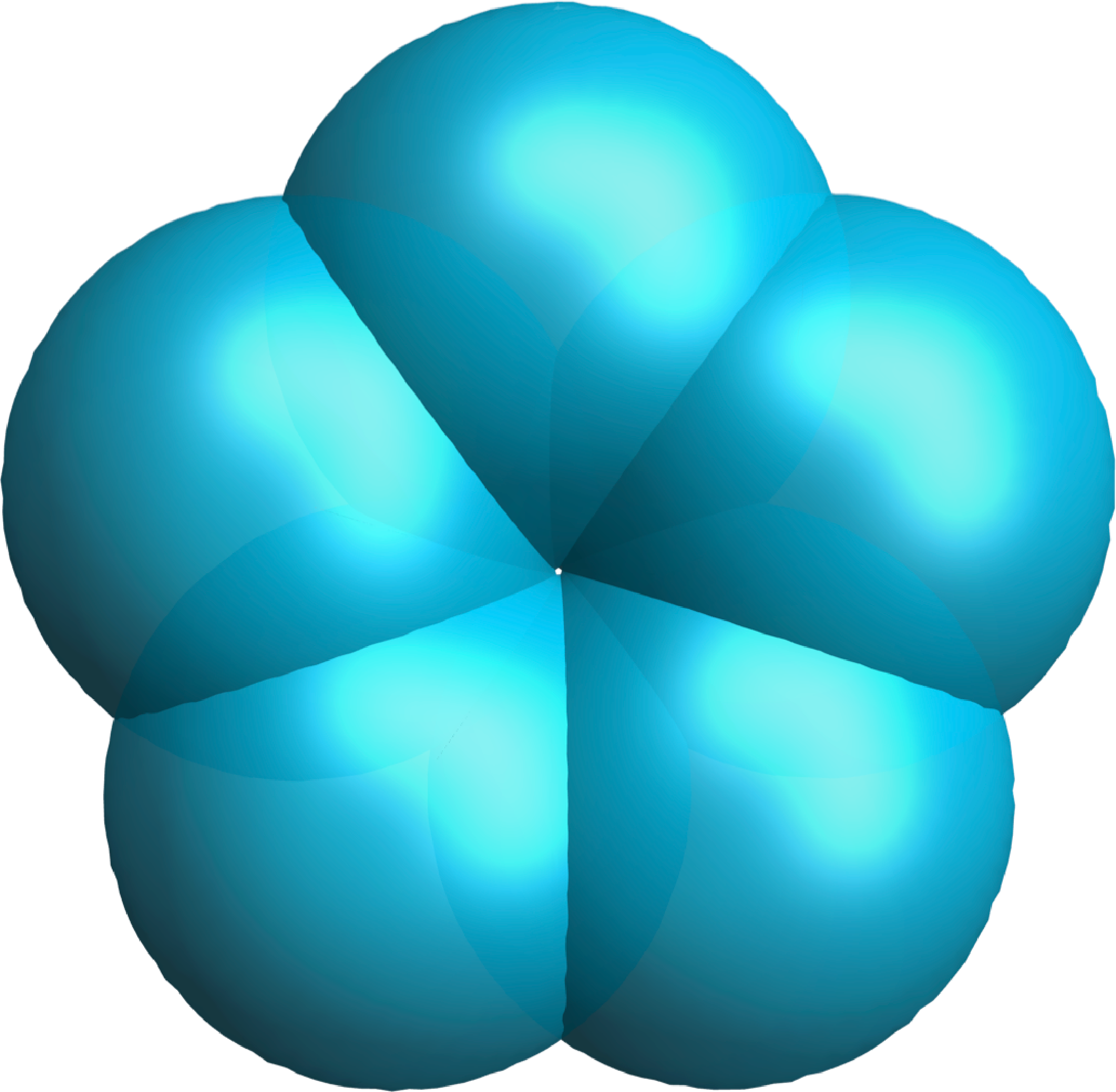}}&
  \resizebox{4.5cm}{!}{\includegraphics{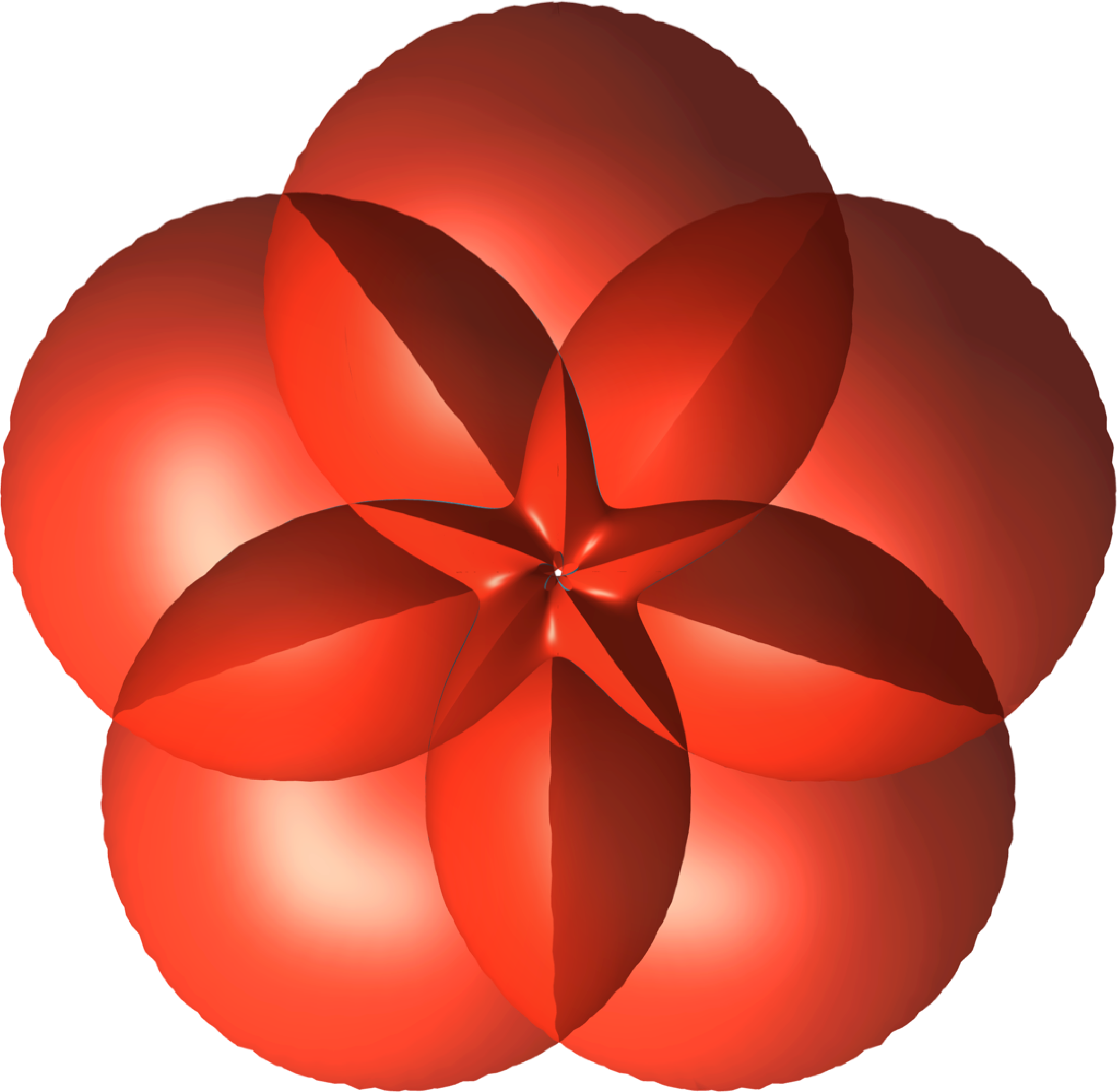}}\\
  {\footnotesize (b1)} & 
  {\footnotesize (b2)} \\ 
 \end{tabular}
\caption{Ruled surfaces in $\S^3$ generated by
K\"{a}hler magnetic geodesics
with charge $\lambda$ in $\LS$, as in
Corollary~\ref{cor:magnetic-ruled-invariant}.
Panel (a1) corresponds to
$\lambda=1/(2\sqrt{2})$ and
$\theta=\frac12\arccos(3/4)$.
Panel (b1) corresponds to
$\lambda=\sqrt{2}/3$ and
$\theta=\frac12\arccos(5/9)$.
Panels (a2) and (b2) show cutaway views of (a1) and (b1), respectively.}
\label{fig:magnetic-ruled-surface}
\end{center}
\end{figure}

\subsection{Total curvature of closed G-projection-regular ruled surfaces}
\label{subsec:ruled-total-curvature}

We conclude by transferring the global inequality of
Section~\ref{sec:fenchel-type} to ruled surfaces.  We first define the
ruling total curvature and characterize simple-rotational ruling families.  
We next prove the sharp lower bound and identify its equality
case.  Finally, we show directly that a nonsingular equality-case ruled
surface is an embedded torus.

Let $F\colon T_L^2\to S^3$ be a closed G-projection-regular ruled surface, where
\[
 T_L^2=(\R/L\mathbb Z)\times(\R/2\pi\mathbb Z),
\]
and let $\gamma\colon \R/L\mathbb Z\to\LS$ be its unit-speed associated curve.

\begin{definition}\label{def:ruling-total-curvature}
The \emph{ruling total curvature} of $F$ is
\begin{equation}\label{eq:ruling-total-curvature}
 \mathcal T(F):=\mathcal T(\gamma).
\end{equation}
Equivalently, it is the product Sasaki length of the canonical lifts of
the right and left ruling Gauss maps.
\end{definition}

The equivalence in the definition follows from
Proposition~\ref{prop:ruling-Gauss-curve}, because the fixed reflection
$\Pi$ and its differential act isometrically on $US^2$.

\begin{definition}\label{def:simple-rotational}
Let $F\colon T_L^2\to S^3$ be a closed G-projection-regular ruled surface
with oriented rulings $\ell_s$.
We call $F$ \emph{simple-rotational} if 
there exists a rotational motion $\{\Phi_s\}$ such that
\[
 \ell_s=\Phi_s(\ell_0)
 \qquad (s\in\R),
\]
and $L$ is the least positive period of the generating motion.
\end{definition}

\begin{proposition}\label{prop:rotational-criterion}
Let $F\colon T_L^2\to S^3$ be a closed G-projection-regular
helicoidal ruled surface. The following are equivalent:
\begin{enumerate}
\item the two ruling Gauss maps are embeddings;
\item the generating helicoidal motion is rotational,
      and $L$ is its least positive period;
\item $F$ is simple-rotational.
\end{enumerate}
\end{proposition}

\begin{proof}
Write the generating motion in the form \eqref{eq:double-rotation}. Since
its associated curve is G-projection-regular, both projections are nonconstant
spherical circles, with angular speeds
$\mu_1=|a+b|$ and $\mu_2=|a-b|$. Closedness gives
\[
 \mu_1L=2\pi m_1,
 \qquad
 \mu_2L=2\pi m_2,
 \qquad m_1,m_2\in\mathbb N,
\]
where each $\gamma_j$ is an $m_j$-fold covering of its image circle.
By \eqref{eq:Gauss-curve-identification}, the corresponding
ruling Gauss map is an embedding if and only if $m_j=1$.

If both ruling Gauss maps are embeddings, then
$|a+b|=|a-b|=2\pi/L$. Squaring gives $4ab=0$, so the generating
motion is rotational. Moreover, its nonzero angular speed is $2\pi/L$,
and therefore its least positive period is $L$.
Conversely, if the generating motion is rotational and $L$ is its
least positive period, then $ab=0$ and
$\mu_1=\mu_2=2\pi/L$. Hence $m_1=m_2=1$, so both
ruling Gauss maps are embeddings.
The equivalence of (2) and (3) follows from
Definition~\ref{def:simple-rotational}.
\end{proof}

\begin{theorem}\label{thm:ruled-Fenchel}
Every closed G-projection-regular ruled surface satisfies
\begin{equation}\label{eq:ruled-Fenchel}
 \mathcal T(F)\geq2\sqrt2\,\pi.
\end{equation}
Equality holds if and only if $F$ is simple-rotational.
Moreover, if
equality holds and $F$ has no singular points, then
$F\colon T_L^2\to S^3$ is an embedding.
\end{theorem}

\begin{proof}
The inequality follows immediately from
Theorem~\ref{prop:fenchel-triplet}. Equality there is equivalent to
the curvature data being constant and both $\gamma_1$ and $\gamma_2$ 
being embeddings.
The first condition is equivalent to
the ruled surface being helicoidal by
Proposition~\ref{prop:helicoidal-characterization}.
The second condition implies that $L$ is the least positive period
of $\gamma$.
By \eqref{eq:Gauss-curve-identification} and
Proposition~\ref{prop:rotational-criterion}, the two conditions
together are equivalent to $F$ being simple-rotational.
The final assertion follows from
Proposition~\ref{prop:equality-embedded} below.
\end{proof}

\begin{proposition}\label{prop:equality-embedded}
Let $F\colon T_L^2\to S^3$ be a closed G-projection-regular ruled surface 
for which equality holds
in \eqref{eq:ruled-Fenchel}.  If $F$ has no singular points, then $F$ is an
embedding.
\end{proposition}

\begin{proof}
Let
$
 \gamma\colon \R/L\Z\to\LS
$
be the associated curve of $F$.
The equality case of Theorem~\ref{prop:fenchel-triplet}
shows that $\gamma$ has constant curvature data and that
both $\gamma_1$ and $\gamma_2$ are embeddings.

By Proposition~\ref{prop:helicoidal-characterization} and
Proposition~\ref{prop:rotational-criterion},
the ruled surface is generated by 
a rotational motion.
After applying a congruence, interchanging the two complex coordinate
planes if necessary, and making a diffeomorphic change of the ruling
parameter, we may write
\begin{equation}\label{eq:rotational-adapted-F}
 F(s,t)=\Phi_s(G(t)),
 \qquad
 \Phi_s(z,w)=(e^{i\omega s}z,w),
 \qquad
 \omega=\frac{2\pi}{L},
\end{equation}
where
$
 S^3=\{(z,w)\in\C^2:|z|^2+|w|^2=1\}
$
and $G$ is an oriented unit-speed geodesic of $S^3$.

Since $G$ is a great circle, there exist orthonormal vectors
$p,q\in\R^4=\C^2$ such that
$G(t)=p\cos t+q\sin t.$
Writing
$p=(z_0,w_0)$ and $q=(z_1,w_1)$,
we obtain
\begin{equation}\label{eq:rotational-profile}
 F(s,t)
 =
 \bigl(e^{i\omega s}z(t),w(t)\bigr)\quad\text{where}\quad
 \left\{ 
  \begin{array}{ll}
    z(t)=z_0\cos t+z_1\sin t,\\
    w(t)=w_0\cos t+w_1\sin t.
  \end{array}
\right.
\end{equation}
We first record two consequences of the nonsingularity of $F$.
Since
$
 F_s
 =
 \bigl(i\omega e^{i\omega s}z(t),0\bigr),
$
we must have
\begin{equation}\label{eq:z-nonzero-rotational}
 z(t)\neq0
 \qquad\text{for all }t.
\end{equation}

We also claim that
\begin{equation}\label{eq:w-prime-nonzero}
 w'(t)\neq0
 \qquad\text{for all }t.
\end{equation}
Indeed, suppose that $w'(t_0)=0$.  
Differentiating
$
 |z(t)|^2+|w(t)|^2=1
$
at $t=t_0$ gives
$
 \operatorname{Re}
 \bigl(\overline{z(t_0)}z'(t_0)\bigr)=0.
$
Since $z(t_0)\neq0$, there exists $\lambda\in\R$ such that
$
 z'(t_0)=i\lambda z(t_0).
$
Consequently,
\[
 F_t(s,t_0)
 =
 \bigl(i\lambda e^{i\omega s}z(t_0),0\bigr)
 =
 \frac{\lambda}{\omega}F_s(s,t_0),
\]
contrary to the nonsingularity of $F$.  This proves
\eqref{eq:w-prime-nonzero}.

Consider the real-linear map
$A\colon\R^2\to\C$
defined by
$A(x,y)=w_0x+w_1y$.
Then
$w(t)=A(\cos t,\sin t)$
and
$ w'(t)=A(-\sin t,\cos t)$.
If $A$ were not injective, then for some $t_0$ we would have
$
 A(-\sin t_0,\cos t_0)=0.
$
This would imply $w'(t_0)=0$.  Hence $A$ is an isomorphism.  It follows
that
$
 w\colon\R/2\pi\Z\to\C
$
is an embedding whose image is a nondegenerate ellipse.

Now suppose that
$
 F(s,t)=F(\widetilde s,\widetilde t).
$
Equality of the second components in
\eqref{eq:rotational-profile} gives
$
 w(t)=w(\widetilde t).
$
Since $w$ is injective,
$
 t=\widetilde t
 \pmod{2\pi}.
$
Using \eqref{eq:z-nonzero-rotational} and comparing the first
components, we obtain
$
 e^{i\omega s}=e^{i\omega\widetilde s}.
$
Since $\omega L=2\pi$, this implies
$s=\widetilde s \pmod{L}.$
Thus $F$ is injective.

By assumption, $F$ has no singular points, and hence it is an
immersion.  Since $T_L^2$ is compact, the injective immersion
$
 F\colon T_L^2\to S^3
$
is an embedding.
\end{proof}

\section*{Declarations}

\subsection*{Funding}
This work was partially supported by JSPS KAKENHI Grants
JP24K06709 (A.H.), JP19K03461 and JP23K03081 (J.I.);
JST SPRING Grant JPMJSP2119 (H.K.);
and the Foundation of Research Fellows,
The Mathematical Society of Japan (H.K.).

\end{document}